\documentclass[11pt,oneside,leqno]{amsproc}
\usepackage[utf8]{inputenc}
\usepackage[russian]{babel}
\usepackage{amssymb,amsmath,amsthm}
\renewcommand{\subsection}{\refstepcounter{subsection}%
\par\bigskip\noindent\textbf{\upshape\thesubsection. }}
\renewcommand{\subsubsection}{\refstepcounter{subsubsection}%
\par\medskip\noindent\textbf{\upshape\thesubsubsection.  }}
\renewcommand{\paragraph}{\refstepcounter{paragraph}%
\par\smallskip\noindent\textbf{\upshape\theparagraph. }}

\numberwithin{equation}{subsection}

\renewcommand{\thesubsection}{\arabic{subsection}}

\newcommand{\Wo}{{\raisebox{0.2ex}{\(\stackrel{\circ}{W}\)}}{}}

\newcommand{\im}{\operatorname{im}}

\title{О чебышёвских свойствах системы собственных функций задачи
Штурма--Лиувилля с сингулярными коэффициентами}
\author{А.~А.~Владимиров\footnote{%
Работа поддержана РФФИ, код проекта~07-01-00283, и фондом INTAS,
код проекта~05-1000008-7883.}}
\begin{document}
\renewcommand{\proofname}{{\upshape Д\,о\,к\,а\,з\,а\,т\,е\,л\,ь\,с\,т\,в\,о.}}
\begin{flushleft}
\normalsize УДК~517.984+517.51
\end{flushleft}
\begin{abstract}
В статье рассматривается сингулярная спектральная задача Штурма--Лиувилля
\begin{gather*}
	-(py')'+(q-\lambda r)y=0,\\
	(U-1)y^{\vee}+i(U+1)y^{\wedge}=0,
\end{gather*}
где функция \(p\in L_{\infty}[0,1]\) равномерно положительна, функция
\(q\in W_2^{-1}[0,1]\) вещественна, весовая функция \(r\in W_2^{-1}[0,1]\)
положительна, а определяющая граничные условия унитарная комплексная матрица
\(U\) размера \(2\times 2\) диагональна. Показывается, что известные
для гладкого случая результаты о чебышёвских свойствах систем собственных
функций остаются справедливыми и в общей ситуации.
\end{abstract}

\maketitle\markboth{}{}
\section{Введение}\label{par:0}
\subsection\label{pt:0:1}
Пусть \(p\) "--- равномерно положительная функция класса \(L_{\infty}[0,1]\),
\(q\) "--- вещественная функция класса \(W_2^{-1}[0,1]\), \(r\) "---
неотрицательная функция класса \(W_2^{-1}[0,1]\) с носителем \([0,1]\),
а \(U\) "--- диагональная унитарная комплексная матрица размера \(2\times 2\).
В настоящей статье будет продолжено начатое в работе \cite{Vl1} исследование
осцилляционных свойств собственных функций граничной задачи
\begin{gather}\label{eq:1}
	-(py')'+(q-\lambda r)y=0,\\ \label{eq:2}
	(U-1)y^{\vee}+i(U+1)y^{\wedge}=0,
\end{gather}
где \(\lambda\in\mathbb C\) "--- спектральный параметр, а векторы
\(y^{\wedge}\) и \(y^{\vee}\) определены в виде
\begin{align*}
	y^{\wedge}&=\begin{pmatrix}y(0)\\ y(1)\end{pmatrix},&
	y^{\vee}&=\begin{pmatrix}y^{[1]}(0)\\ -y^{[1]}(1)\end{pmatrix}.
\end{align*}

Ввиду негладкости коэффициентов дифференциального выражения в левой части
уравнения~\eqref{eq:1}, постановка задачи~\eqref{eq:1}, \eqref{eq:2}
нуждается в уточнении. Проводя его, мы, как и в работе \cite{Vl1}, будем
исходить из аппроксимативной точки зрения (см. \cite{SaSh,Vl}).

Обозначим через \(\mathfrak H_U\) гильбертового пространство
\[
	\{y\in W_2^1[0,1]\mid y^{\wedge}\in\im(U-1)\},
\]
снабжённое обычной нормой
\[
	(\forall y\in\mathfrak H_U)\qquad \|y\|^2_{\mathfrak H_U}=
	\int\limits_0^1\left\{|y'|^2+|y|^2\right\}\,dx.
\]
Рассмотрим оператор вложения \(I:\mathfrak H_U\to L_2[0,1]\) и обозначим
через \(\mathfrak H'_U\) пополнение пространства \(L_2[0,1]\) по норме
\(\|y\|_{\mathfrak H'_U}\rightleftharpoons \|I^*y\|_{\mathfrak H_U}\).
Непосредственно из этого определения вытекает возможность непрерывного
продолжения оператора \(I^*\) до изометрии \(I^+:\mathfrak H'_U\to
\mathfrak H_U\). Граничную задачу~\eqref{eq:1}, \eqref{eq:2} мы будем теперь
понимать как задачу о спектре линейного операторного пучка
\(T:\mathfrak H_U\to\mathfrak H_U'\), имеющего вид
\[
	(\forall\lambda\in\mathbb C)\:(\forall y\in\mathfrak H_U)\qquad
	\langle I^+T(\lambda)y,y\rangle_{\mathfrak H_U}=
	\int\limits_0^1 p\,|y'|^2\,dx+
	\int\limits_0^1 (q-\lambda r)\cdot |y|^2\,dx+
	\langle Vy^{\wedge}, y^{\wedge}\rangle_{\mathbb C^2},
\]
где \(V\) "--- диагональная эрмитова матрица размера \(2\times 2\)
с элементами
\[
	(\forall k\in\{1,2\})\qquad V_{kk}=\left\{\begin{aligned}
	-&\ctg\dfrac{\arg U_{kk}}{2},&&U_{kk}\neq 1,\\
	&0,&&U_{kk}=1.\end{aligned}\right.
\]

\subsection
В работе \cite{Vl1} было показано, что спектр задачи
\ref{pt:0:1}\,\eqref{eq:1}, \ref{pt:0:1}\,\eqref{eq:2} может быть представлен
в виде возрастающей последовательности простых собственных значений
\[
	\lambda_1<\lambda_2<\ldots<\lambda_n<\ldots
\]
Было показано также, что отвечающие этим собственным значениям собственные
функции \(y_n\in\mathfrak H_U\) имеют по \(n-1\)~нулей на интервале \((0,1)\),
причём нули соседних собственных функций перемежаются. Между тем, для случая
гладкости коэффициентов дифференциального уравнения
\ref{pt:0:1}\,\eqref{eq:1} имеются и существенно более сильные результаты,
связанные с теорией осцилляционных матриц и ядер (см., например, \cite{GK},
\cite{BP}). В частности, известно, что произвольная ненулевая линейная
комбинация
\[
	\sum\limits_{k=n}^N\alpha_k\,y_k
\]
вещественных собственных функций \(y_k\) с вещественными коэффициентами
\(\alpha_k\) имеет на интервале \((0,1)\) не менее \(n-1\) перемен знака
и не более \(N-1\) нулей. Цель настоящей статьи заключается в установлении
справедливости указанного факта и в общей ситуации.

\subsection\label{pt:0:3}
Заметим, что в рассматриваемом нами случае диагональности матрицы \(U\)
пространство \(\mathfrak H_U\) может иметь только один из четырёх следующих
видов:
\begin{gather*}
	\Wo_2^1[0,1],\\
	\{y\in W_2^1[0,1]\mid y(1)=0\},\\
	\{y\in W_2^1[0,1]\mid y(0)=0\},\\
	W_2^1[0,1].
\end{gather*}
Это означает, что за счёт возмущений потенциала \(q\) сосредоточенными
в точках \(0\) и \(1\) дельта-функциями любую задачу рассматриваемого нами
типа можно переписать в виде аналогичной задачи с граничными условиями одного
из следующих четырёх видов:
\begin{gather}\label{eq:0:1}
	y(0)=y(1)=0,\\ \label{eq:0:2}
	y^{[1]}(0)=y(1)=0,\\ \label{eq:0:3}
	y(0)=y^{[1]}(1)=0,\\ \label{eq:0:4}
	y^{[1]}(0)=y^{[1]}(1)=0.
\end{gather}
Заметим также, что рассмотрение задачи с граничными условиями~\eqref{eq:0:3}
может быть заменой переменной \(t\mapsto \tau\rightleftharpoons 1-t\) сведено
к рассмотрению задачи с граничными условиями~\eqref{eq:0:2}. Поэтому
в дальнейшем мы, как правило, будем рассматривать только задачи с граничными
условиями \eqref{eq:0:1}, \eqref{eq:0:2} и \eqref{eq:0:4}.

\subsection
Структура оставшейся части статьи такова. В параграфе~\ref{par:1} нами
устанавливаются вспомогательные результаты о возможности исключения потенциала
\(q\). В параграфе~\ref{par:2} устанавливается знакорегулярность резольвенты
задачи \ref{pt:0:1}\,\eqref{eq:1}, \ref{pt:0:1}\,\eqref{eq:2}. Наконец,
в параграфе~\ref{par:3} проводится доказательство основных результатов статьи.

При ссылках на разделы статьи, не принадлежащие параграфу, внутри которого
даётся ссылка, дополнительно указывается номер параграфа. При ссылках
на формулы, не принадлежащие пункту, внутри которого даётся ссылка,
дополнительно указывается номер пункта.


\section{Исключение потенциала}\label{par:1}
\subsection
Хорошо известно (см, например, \cite[\S\,14]{KK}), что изучение задачи
\ref{par:0}.\ref{pt:0:1}\,\eqref{eq:1}, \ref{par:0}.\ref{pt:0:1}\,\eqref{eq:2}
с гладкими коэффициентами \(p\), \(q\) и \(r\) может быть посредством замены
переменной сведено к изучению аналогичной задачи с потенциалом \(q\equiv 0\).
В настоящем параграфе соответствующая техника будет распространена на случай
коэффициентов рассматриваемого нами более общего вида.

На протяжении настоящего параграфа мы будем предполагать зафиксированным
число \(\xi\in\mathbb R\), для которого оператор \(I^+T(\xi)\) является
положительным. Существование такого числа установлено, например, в работе
\cite[\S\,3]{Vl1}.

\subsection
Заметим (см. \cite{SaSh, Vl}), что существуют функция \(\omega\in L_2[0,1]\)
и число \(\omega_1\in\mathbb R\), удовлетворяющие тождеству
\begin{equation}\label{eq:1:0}
	(\forall y\in W_2^1[0,1])\qquad \int\limits_0^1 (q-\xi r)\cdot
	\overline{y}\,dx=-\int\limits_0^1 \omega\cdot\overline{y'}\,dx+
	\omega_1\cdot\overline{y(1)}.
\end{equation}
При этом имеют место следующие два факта:

\subsubsection\label{prop:1:1}
{\itshape Существует вектор-функция \(Y\in W_2^1[0,1]\times W_1^1[0,1]\),
обладающая равномерно на отрезке \([0,1]\) положительной первой компонентой
\(Y_1\) и удовлетворяющая уравнениям
\begin{gather}\label{eq:1:1}
	\dfrac{dY}{dt}=\begin{pmatrix}\omega/p& 1/p\\
	-\omega^2/p& -\omega/p\end{pmatrix}\cdot Y,\\ \notag
	\int\limits_0^1\dfrac{dx}{Y_1^2}=1.
\end{gather}
При этом в случае, когда граничные условия
\ref{par:0}.\ref{pt:0:1}\,\eqref{eq:2} имеют один из видов
\ref{par:0}.\ref{pt:0:3}\,\eqref{eq:0:2} или
\ref{par:0}.\ref{pt:0:3}\,\eqref{eq:0:4}, можно дополнительно потребовать
выполнения равенства
\begin{equation}\label{eq:1:2}
	Y_2(0)=0.
\end{equation}
}

\begin{proof}
Здесь мы будем опираться на установленные в работе \cite[\S\,2]{Vl1}
результаты о точках, сопряжённых точкам \(0\) и \(1\) относительно оператора
\(T(\xi)\).

Если граничные условия имеют один из видов
\ref{par:0}.\ref{pt:0:3}\,\eqref{eq:0:2} или
\ref{par:0}.\ref{pt:0:3}\,\eqref{eq:0:4}, то для доказательства
рассматриваемого утверждения достаточно заметить, что полуинтервал \((0,1]\)
не содержит точек, сопряжённых точке \(0\) относительно оператора \(T(\xi)\).

Если граничные условия имеют вид \ref{par:0}.\ref{pt:0:3}\,\eqref{eq:0:1},
то зафиксируем решение \(Z\in W_2^1[0,1]\times W_1^1[0,1]\) отвечающей
уравнению \eqref{eq:1:1} начальной задачи
\[
	Y_1(1)=Y_2(1)+1=0.
\]
Поскольку полуинтервал \([0,1)\) не содержит точек, сопряжённых точке \(1\)
относительно оператора \(T(\xi)\), то функция \(Z_1\) не может иметь нулей
на этом полуинтервале. Отсюда, в свою очередь, следует, что для любого
значения \(C\geqslant 0\) квадратичная форма
\[
	\int\limits_0^1 p\,|y'|^2\,dx+\int\limits_0^1 (q-\xi r)\cdot
	|y|^2\,dx+\dfrac{Z_2(0)+C}{Z_1(0)}\cdot |y(0)|^2
\]
является неотрицательной на пространстве
\[
	\{y\in W_2^1[0,1]\mid y(1)=0\}.
\]
Последнее означает, что для любого значения \(C\geqslant 0\) первая
компонента решения отвечающей уравнению \eqref{eq:1:1} начальной задачи
\[
	Y_1(0)-Z_1(0)=Y_2(0)-Z_2(0)-C=0
\]
является положительной на интервале \((0,1)\). Для завершения доказательства
теперь достаточно заметить, что при \(C>0\) теорема единственности
\cite[\S\,16, Теорема~1]{Na} гарантирует выполнение неравенства
\(Y_1(1)\neq 0\).
\end{proof}

\subsubsection\label{prop:1:2}
{\itshape Пусть вектор-функция \(Y\in W_2^1[0,1]\times W_1^1[0,1]\) обладает
указанными в утверждении \ref{prop:1:1} свойствами. Тогда выполняется
тождество
\begin{equation}\label{eq:1:3}
	(\forall y\in\mathfrak H_U)\qquad
	\int\limits_0^1 p\,Y_1'\overline{y'}\,dx+
	\int\limits_0^1 (q-\xi r)\cdot (Y_1\overline{y})\,dx=
	[Y_2(1)+\omega_1\cdot Y_1(1)]\cdot\overline{y(1)}.
\end{equation}
При этом в случае, когда граничные условия
\ref{par:0}.\ref{pt:0:1}\,\eqref{eq:2} имеют вид
\ref{par:0}.\ref{pt:0:3}\,\eqref{eq:0:4}, выполняется также неравенство
\begin{equation}\label{eq:1:4}
	\dfrac{Y_2(1)}{Y_1(1)}+\omega_1>0.
\end{equation}
}

\begin{proof}
Из соотношений \eqref{eq:1:0} и \eqref{eq:1:1} следует, что при любом
\(y\in\mathfrak H_U\) выполняются равенства
\begin{flalign*}
	&&\int\limits_0^1 pY_1'\overline{y'}\,dx+
	\int\limits_0^1 (q-\xi r)\cdot Y_1\overline{y}\,dx&=
	\int\limits_0^1 [\omega Y_1+Y_2]\,\overline{y'}\,dx-
	\int\limits_0^1 \omega\cdot[Y_1'\overline{y}+
	Y_1\overline{y'}]\,dx+\omega_1 Y_1(1)\overline{y(1)}&&\\
	&& &=\int\limits_0^1 [Y_2\overline{y'}+Y_2'\overline{y}]\,dx+
	\omega_1 Y_1(1)\overline{y(1)}\\
	&& &=[Y_2(1)+\omega_1\cdot Y_1(1)]\cdot \overline{y(1)}-
	Y_2(0)\overline{y(0)}.
\end{flalign*}
Отсюда и из предполагаемых выполненными для случая граничных условий
\ref{par:0}.\ref{pt:0:3}\,\eqref{eq:0:2} или
\ref{par:0}.\ref{pt:0:3}\,\eqref{eq:0:4} равенств \eqref{eq:1:2} немедленно
вытекает справедливость тождества \eqref{eq:1:3}. Неравенство \eqref{eq:1:4}
представляет собой несложное следствие тождества \eqref{eq:1:3} и факта
положительности оператора \(I^+T(\xi)\).
\end{proof}

\subsection\label{pt:1:3}
Зафиксируем вектор-функцию \(Y\in W_2^1[0,1]\times W_1^1[0,1]\) с описанными
в утверждении~\ref{prop:1:1} свойствами и сопоставим ей возрастающую функцию
\(\tau:[0,1]\to [0,1]\) вида
\[
	(\forall t\in [0,1])\qquad \tau(t)=\int\limits_0^t\dfrac{dx}{Y_1^2}.
\]
С функцией \(\tau\) мы в дальнейшем будем связывать ограниченно обратимый
оператор замены переменной \(J:\mathfrak H_U\to \mathfrak H_U\), имеющий
вид
\[
	(\forall y\in\mathfrak H_U)\:(\forall t\in [0,1])\qquad
	[Jy](\tau(t))=y(t),
\]
а также получаемые продолжением этого оператора по непрерывности операторы
\(J_{\infty}:L_{\infty}[0,1]\to L_{\infty}[0,1]\) и \(J_{-1}:W_2^{-1}[0,1]\to
W_2^{-1}[0,1]\). Кроме того, мы будем рассматривать ограниченно обратимый
оператор \(S:\mathfrak H_U\to\mathfrak H_U\) вида
\begin{equation}\label{eq:1:5}
	(\forall y\in\mathfrak H_U)\qquad Sy=J(y/Y_1).
\end{equation}
Имеет место следующий факт:

\subsubsection\label{prop:1:3}
{\itshape Пусть \(\hat T:\mathfrak H_U\to\mathfrak H'_U\) "--- линейный
операторный пучок, отвечающий дифференциальному уравнению
\[
	-(\hat py')'-(\lambda-\xi)\hat ry=0,
\]
где положено \(\hat p\rightleftharpoons J_{\infty}p\)
и \(\hat r\rightleftharpoons J_{-1}(Y_1^4\,r)\). Пусть также в случае, когда
исходный пучок \(T\) определяется граничными условиями
\ref{par:0}.\ref{pt:0:3}\,\eqref{eq:0:1} или
\ref{par:0}.\ref{pt:0:3}\,\eqref{eq:0:2}, пучок \(\hat T\) определяется теми же
граничными условиями, а в случае, когда исходный пучок \(T\) определяется
граничными условиями \ref{par:0}.\ref{pt:0:3}\,\eqref{eq:0:4}, пучок \(\hat T\)
определяется граничными условиями
\[
	y^{[1]}(0)=y^{[1]}(1)+\left[\dfrac{Y_2(1)}{Y_1(1)}+\omega_1\right]\cdot
	y(1)=0.
\]
Тогда выполняется тождество
\[
	(\forall\lambda\in\mathbb C)\qquad S^*I^+\hat T(\lambda)S=I^+T(\lambda).
\]
}

\begin{proof}
Заметим, что непосредственно из определения оператора \(S\) вытекает
справедливость тождества
\[
	(\forall y\in\mathfrak H_U)\:(\forall t\in [0,1])\qquad
	[Sy]'(\tau(t))\cdot\dfrac{1}{Y_1^2(t)}=
	\dfrac{Y_1(t)\cdot y'(t)-Y_1'(t)\cdot y(t)}{Y_1^2(t)}.
\]
Поэтому для любой функции \(y\in\mathfrak H_U\) выполняются равенства
\begin{flalign*}
	&& \int\limits_0^1\hat p\,\bigl|[Sy]'\bigr|^2\,dx&=
	\int\limits_0^1 p\,|Y_1y'-Y_1'y|^2\,\dfrac{dx}{Y_1^2}&&\\
	&& &=\int\limits_0^1 p\,|y'|^2\,dx-\int\limits_0^1 pY_1'
	\overline{\left(\dfrac{|y|^2}{Y_1}\right)'}\,dx.
\end{flalign*}
Аналогично, для любой функции \(y\in\mathfrak H_U\) выполняются равенства
\begin{flalign*}
	&& \int\limits_0^1\hat r\cdot \bigl|[Sy]\bigr|^2\,dx&=
	\int\limits_0^1 (Y_1^4\,r)\cdot \left|\dfrac{y}{Y_1}\right|^2\,
	\dfrac{dx}{Y_1^2}&&\\
	&& &=\int\limits_0^1 r\cdot |y|^2\,dx.
\end{flalign*}
Объединяя эти факты со следующим из утверждения \ref{prop:1:2} тождеством
\[
	(\forall y\in\mathfrak H_U)\qquad -\int\limits_0^1 pY_1'
	\overline{\left(\dfrac{|y|^2}{Y_1}\right)'}\,dx=\int\limits_0^1
	(q-\xi r)\cdot |y|^2\,dx-\left[\dfrac{Y_2(1)}{Y_1(1)}+
	\omega_1\right]\cdot |y(1)|^2,
\]
убеждаемся в справедливости тождества
\[
	(\forall\lambda\in\mathbb C)\:(\forall y\in\mathfrak H_U)\qquad
	\langle I^+\hat T(\lambda)\,[Sy],[Sy]\rangle_{\mathfrak H_U}=
	\langle I^+T(\lambda)y,y\rangle_{\mathfrak H_U},
\]
равносильного доказываемому утверждению.
\end{proof}

\subsection
Отметим, что для частного случая, выделяемого условиями \(p\in BV[0,1]\),
\(q\geqslant 0\) и \(U=1\), утверждение о возможности исключения потенциала
было сформулировано в недавней работе \cite{Pok}.


\section{Знакорегулярность резольвенты}\label{par:2}
\subsection\label{pt:2:1}
В дальнейшем через \(C[0,1]\) мы будем обозначать пространство равномерно
непрерывных на отрезке \([0,1]\) вещественных функций. При этом, как обычно,
мы будем говорить, что функция \(f\in C[0,1]\) \emph{имеет не менее \(n\)
перемен знака на интервале \((0,1)\)}, если найдётся набор
\(\{x_k\}_{k=1}^{n+1}\) из \(n+1\)~точек
\begin{equation}\label{eq:3:0}
	0<x_1<x_2<\ldots<x_{n+1}<1,
\end{equation}
удовлетворяющий условию
\[
	(\forall k\in\{1,\ldots, n\})\qquad f(x_k)\cdot f(x_{k+1})<0.
\]
Аналогично, мы будем говорить, что функция \(f\in C[0,1]\) \emph{имеет
не менее \(n\) псевдонулей на интервале \((0,1)\)}, если найдутся вещественное
число \(\varepsilon>0\) и набор \(\{x_k\}_{k=1}^{n+1}\) из \(n+1\)~точек вида
\eqref{eq:3:0}, удовлетворяющие условиям
\begin{gather*}
	(\forall k\in\{1,\ldots,n+1\})\qquad |f(x_k)|>\varepsilon,\\
	(\forall k\in\{1,\ldots,n\})\,(\exists t\in (x_k,x_{k+1}))\qquad
	|f(t)|<\varepsilon.
\end{gather*}
Очевидна справедливость следующих двух утверждений:

\subsubsection\label{prop:2:0}
{\itshape Пусть \(n\) "--- натуральное число, и пусть функция \(f\in C[0,1]\)
имеет не менее \(n\) перемен знака на интервале \((0,1)\). Тогда функция \(f\)
имеет не менее \(n\) псевдонулей на интервале \((0,1)\).
}

\subsubsection\label{prop:2:1}
{\itshape Пусть \(n\) "--- натуральное число. Тогда множество функций
\(f\in C[0,1]\), имеющих не менее \(n\) псевдонулей на интервале \((0,1)\),
является открытым в топологии пространства \(C[0,1]\).
}

\subsection
Вещественный оператор \(K:\mathfrak H_U\to\mathfrak H_U\) мы в дальнейшем
будем называть \emph{усиленно знакорегулярным}, если для любых натурального
числа \(n\) и вещественной функции \(y\in\mathfrak H_U\), для которых функция
\(Ky\in\mathfrak H_U\) имеет не менее \(n\) псевдонулей на интервале
\((0,1)\), функция \(y\) имеет не менее \(n\) перемен знака на этом интервале.
Введённое понятие примыкает к известному понятию знакорегулярного оператора
(см. \cite{LS}, \cite[\S\,4]{BP}).

Очевидным следствием утверждения~\ref{prop:2:0} является такое утверждение:

\subsubsection\label{prop:2:2}
{\itshape Пусть \(K_1:\mathfrak H_U\to\mathfrak H_U\)
и \(K_2:\mathfrak H_U\to\mathfrak H_U\) "--- два усиленно знакорегулярных
оператора. Тогда оператор \(K_1K_2\) также является усиленно знакорегулярным.
}

\bigskip
Кроме того, имеет место следующий факт:

\subsubsection\label{prop:2:3}
{\itshape Пусть оператор \(K:\mathfrak H_U\to\mathfrak H_U\) является пределом
некоторой последовательности \(\{K_k\}_{k=1}^{\infty}\) усиленно знакорегулярных
операторов в смысле сильной операторной топологии. Тогда он также является
усиленно знакорегулярным.
}

\begin{proof}
Зафиксируем произвольные натуральное число \(n\) и вещественную функцию
\(y\in\mathfrak H_U\), для которых функция \(Ky\in\mathfrak H_U\) имеет
не менее \(n\) псевдонулей на интервале \((0,1)\). Из утверждения
\ref{prop:2:1} и непрерывного характера вложения вещественной части
пространства \(\mathfrak H_U\) в пространство \(C[0,1]\) следует, что
найдётся номер \(m\geqslant 1\), для которого функция \(K_my\) также будет
иметь не менее \(n\) псевдонулей на интервале \((0,1)\). В таком случае
функция \(y\) имеет не менее \(n\) перемен знака на интервале \((0,1)\).
\end{proof}

\subsection\label{pt:2:2}
Основным результатом настоящего параграфа является следующий факт:

\subsubsection\label{prop:2:4}
{\itshape Пусть определяющее пучок \(T\) дифференциальное уравнение
\ref{par:0}.\ref{pt:0:1}\,\eqref{eq:1} удовлетворяет условию \(q=0\),
а граничные условия \ref{par:0}.\ref{pt:0:1}\,\eqref{eq:2} имеют один из видов
\ref{par:0}.\ref{pt:0:3}\,\eqref{eq:0:1},
\ref{par:0}.\ref{pt:0:3}\,\eqref{eq:0:2} или
\begin{equation}\label{eq:3:1}
	y^{[1]}(0)=y^{[1]}(1)+C\cdot y(1)=0,\qquad \text{где } C>0.
\end{equation}
Тогда оператор \(R:\mathfrak H_U\to\mathfrak H_U\) вида
\begin{equation}\label{eq:3:2}
	R\rightleftharpoons [T(0)]^{-1}\cdot\dfrac{dT}{d\lambda}
\end{equation}
является усиленно знакорегулярным.
}

\begin{proof}
Из утверждения \ref{prop:2:3} следует, что достаточно рассмотреть случай,
когда коэффициенты \(p\) и \(r\) дифференциального уравнения
\ref{par:0}.\ref{pt:0:1}\,\eqref{eq:1} являются равномерно положительными
гладкими функциями. Именно такой случай мы и будем изучать в оставшейся
части доказательства. Заметим также, что для доказательства рассматриваемого
утверждения достаточно установить, что для любых натурального числа \(n\)
и вещественной функции \(y\in\mathfrak H_U\), для которых функция
\(u\rightleftharpoons Ry\in\mathfrak H_U\) имеет не менее \(n\) псевдонулей
на интервале \((0,1)\), функция \(y\) имеет не менее \(n\) перемен знака
на этом интервале. Поэтому в оставшейся части доказательства натуральное число
\(n\) и вещественную функцию \(y\in\mathfrak H_U\) мы будем предполагать
зафиксированными.

Пусть граничные условия имеют вид \ref{par:0}.\ref{pt:0:3}\,\eqref{eq:0:1}.
В таком случае из теоремы Лагранжа о конечном приращении легко выводится,
что функция \(pu'\in C^1[0,1]\) имеет не менее \(n+1\) перемен знака
на интервале \((0,1)\). Применяя эту теорему ещё раз, устанавливаем,
что функция \(y=(pu')'/r\) имеет не менее \(n\) перемен знака на том же
интервале. Тем самым, для случая граничных условий
\ref{par:0}.\ref{pt:0:3}\,\eqref{eq:0:1} доказываемое утверждение справедливо.

Пусть граничные условия имеют вид \ref{par:0}.\ref{pt:0:3}\,\eqref{eq:0:2}.
В этом случае из теоремы Лагранжа выводится, что функция \(pu'\in C^1[0,1]\)
имеет не менее \(n\) перемен знака на интервале \((0,1)\). Отсюда, с учётом
равенства \([pu'](0)=0\), следует, что функция \(y=(pu')'/r\) также имеет
не менее \(n\) перемен знака на интервале \((0,1)\). Тем самым, для случая
граничных условий \ref{par:0}.\ref{pt:0:3}\,\eqref{eq:0:2} доказываемое
утверждение также справедливо.

Пусть теперь граничные условия имеют вид \eqref{eq:3:1}. Из положительности
коэффициента \(C\) следует, что для любых двух вещественных чисел
\(\varepsilon>0\) и \(x\in (0,1)\), удовлетворяющих условиям
\begin{gather*}
	|u(x)|<\varepsilon,\\
	(\exists t\in (x,1))\qquad |u(t)|>\varepsilon,
\end{gather*}
найдутся две точки \(t_{\pm}\in (x,1)\), удовлетворяющие неравенствам
\[
	\pm[pu'](t_{\pm})>0.
\]
Этот факт с очевидностью означает, что функция \(pu'\in C^1[0,1]\) имеет
не менее \(n\) перемен знака на интервале \((0,1)\). Отсюда, с учётом
равенства \([pu'](0)=0\), следует, что функция \(y=(pu')'/r\) также имеет
не менее \(n\) перемен знака на этом интервале. Тем самым, справедливость
доказываемого утверждения установлена полностью.
\end{proof}


\section{Чебышёвские свойства систем собственных функций}\label{par:3}
\subsection
Зафиксируем последовательность \(\{y_k\}_{k=1}^{\infty}\) вещественных
собственных функций задачи \ref{par:0}.\ref{pt:0:1}\,\eqref{eq:1},
\ref{par:0}.\ref{pt:0:1}\,\eqref{eq:2}, отвечающих расположенным в порядке
возрастания собственным значениям
\[
	\lambda_1<\lambda_2<\ldots<\lambda_n<\ldots
\]
Имеют место следующие два факта:

\subsubsection\label{prop:3:1}
{\itshape Пусть \(n\) и \(N\) "--- два натуральных числа, удовлетворяющих
условию \(1\leqslant n\leqslant N\). Пусть также \(y\in\mathfrak H_U\) "---
функция, допускающая представление в виде
\[
	y=\sum\limits_{k=n}^N\alpha_k\,y_k,
\]
где \(\{\alpha_k\}_{k=n}^N\) "--- набор вещественных чисел, удовлетворяющий
условию \(\alpha_n\neq 0\). Тогда функция \(y\) имеет на интервале \((0,1)\)
не менее \(n-1\) перемен знака.
}

\begin{proof}
Из утверждений \ref{par:1}.\ref{prop:1:3}, \ref{par:1}.\ref{prop:1:2} и
характера определённого соотношением \ref{par:1}.\ref{pt:1:3}\,\eqref{eq:1:5}
преобразования \(S\) следует, что достаточно рассмотреть случай, когда
выполняется равенство \(q=0\), а определяющие пучок \(T\) граничные условия
имеют один из видов \ref{par:0}.\ref{pt:0:3}\,\eqref{eq:0:1},
\ref{par:0}.\ref{pt:0:3}\,\eqref{eq:0:2} или
\ref{par:2}.\ref{pt:2:2}\,\eqref{eq:3:1}. Именно такой случай мы и будем
изучать в оставшейся части доказательства.

Рассмотрим определённый соотношением \ref{par:2}.\ref{pt:2:2}\,\eqref{eq:3:2}
оператор \(R\). Ввиду простоты собственных значений пучка \(T\),
последовательность \(\{\lambda_n^m\,R^my\}_{m=1}^{\infty}\) сходится
к собственной функции \(\alpha_n\,y_n\). Согласно результатам работы
\cite{Vl1}, такая собственная функция имеет не менее \(n-1\) псевдонулей
на интервале \((0,1)\). Справедливость доказываемого утверждения вытекает
теперь из утверждений \ref{par:2}.\ref{prop:2:1}, \ref{par:2}.\ref{prop:2:2},
\ref{par:2}.\ref{prop:2:4} и непрерывного характера вложения вещественной
части пространства \(\mathfrak H_U\) в пространство \(C[0,1]\).
\end{proof}

\subsubsection\label{prop:3:2}
{\itshape Пусть \(N\) "--- натуральное число, удовлетворяющее условию
\(N\geqslant 1\). Пусть также \(y\in\mathfrak H_U\) "--- функция, допускающая
представление в виде
\[
	y=\sum\limits_{k=1}^N\alpha_k\,y_k,
\]
где \(\{\alpha_k\}_{k=1}^N\) "--- набор вещественных чисел, удовлетворяющий
условию \(\alpha_N\neq 0\). Тогда функция \(y\) имеет на интервале \((0,1)\)
не более \(N-1\) различных нулей.
}

\begin{proof}
Здесь, как и при доказательстве утверждения \ref{prop:3:1}, достаточно
ограничиться рассмотрением случая, когда выполняется равенство \(q=0\),
а определяющие пучок \(T\) граничные условия имеют один из видов
\ref{par:0}.\ref{pt:0:3}\,\eqref{eq:0:1},
\ref{par:0}.\ref{pt:0:3}\,\eqref{eq:0:2} или
\ref{par:2}.\ref{pt:2:2}\,\eqref{eq:3:1}.

Рассмотрим набор \(\{u_m\}_{m=0}^{N}\) принадлежащих пространству
\(\mathfrak H_U\) функций вида
\[
	u_m=\sum\limits_{k=1}^N (\lambda_k^m\alpha_k)\,y_k.
\]
Функции из этого набора с очевидностью удовлетворяют условию
\[
	(\forall m\in\{1,\ldots,N\})\qquad u_{m-1}=Ru_m,
\]
где \(R\) "--- оператор, определённый соотношением
\ref{par:2}.\ref{pt:2:2}\,\eqref{eq:3:2}. Предположим, что для некоторого
номера \(m\in \{1,2,\ldots,N\}\) функция \(u_{m-1}\) тождественно обращается
в нуль на некотором интервале \(\Delta\subset (0,1)\). Тогда для любой функции
\(v\in\Wo_2^1[0,1]\), удовлетворяющей условию \(\operatorname{supp} v\subset
\Delta\), выполняются равенства
\begin{flalign*}
	&&\int\limits_0^1 r\,u_m\overline{v}\,dx&=
	\int\limits_0^1 p\,u_{m-1}'\overline{v'}\,dx&&\\
	&& &=0.
\end{flalign*}
Ввиду максимальности носителя функции \(r\) это означает, что функция
\(u_m\) также тождественно обращается в нуль на интервале \(\Delta\).
Тем самым, исчезновение функции \(y=u_0\) на некотором интервале
\(\Delta\subset (0,1)\) повлекло бы за собой исчезновение на этом же интервале
всех функций из набора \(\{u_k\}_{k=0}^{N}\), а потому и собственной функции
\(y_N\). Однако последнее заведомо невозможно (см., например,
\cite[\S\,2]{Vl1}). Следовательно, функция \(y\) могла бы иметь \(N\)
различных нулей на интервале \((0,1)\) лишь в том случае, если бы она имела
не менее \(N\) псевдонулей на этом интервале. Из утверждения
\ref{par:2}.\ref{prop:2:4} теперь следует, что для завершения доказательства
рассматриваемого утверждения достаточно установить, что функция \(u_1\) не может
иметь более чем \(N-1\) перемен знака на интервале \((0,1)\).

Заметим, что в случае гладкости коэффициентов \(p\gg 0\), \(q\equiv 0\)
и \(r\gg 0\) никакая вещественная линейная комбинация первых \(N\)
собственных функций задачи \ref{par:0}.\ref{pt:0:1}\,\eqref{eq:1},
\ref{par:0}.\ref{pt:0:1}\,\eqref{eq:2} заведомо не может иметь на интервале
\((0,1)\) более чем \(N-1\) перемен знака (см., например,
\cite[Глава~\mbox{IV}, \S\,10]{GK}). Однако произвольная задача рассматриваемого
нами типа может быть в смысле сильной операторной топологии над пространством
\(\mathfrak H_U\) приближена аналогичными задачами с гладкими коэффициентами.
Поэтому сказанное означает невозможность для функции \(u_1\) иметь на интервале
\((0,1)\) более чем \(N-1\) перемен знака и в общем случае. Тем самым,
доказываемое утверждение справедливо.
\end{proof}

\subsection
Отметим, что в действительности свойство максимальности носителя весовой
функции \(r\) не является принципиальным для большинства установленных ранее
результатов. А именно, утверждения \ref{par:1}.\ref{prop:1:3},
\ref{par:2}.\ref{prop:2:4} и \ref{prop:3:1} с очевидностью остаются
справедливыми и для случая, когда носитель весовой функции \(r\) произволен,
а оператор \(I^+T(\xi)\) является положительным при некотором значении
\(\xi\in\mathbb R\). Что касается утверждения \ref{prop:3:2}, то оно
в указанном случае также остаётся справедливым, если заменить в его
формулировке утверждение о несуществовании \(N\) различных нулей функции
\(y\) на интервале \((0,1)\) более слабым утверждением о несуществовании
\(N\) различных связных компонент множества нулей этой функции, не содержащих
точек \(0\) и \(1\).

\end{document}